\documentclass[11pt]{article}

\usepackage{lscape,tikz}

\definecolor{darkblue}{rgb}{0.2,0.2,0.71}
\usepackage{todonotes}
\usepackage{float}
\usepackage{colortbl}
\usepackage{multicol}
\usepackage{wrapfig}
\usepackage{color}

\usepackage{epsfig,cancel,ulem, amsthm,amssymb,amsmath}
\usepackage{color,soul}
\usepackage{graphicx}
\usepackage{yhmath}
\usepackage{mathdots}
\usepackage{MnSymbol}

\usetikzlibrary{decorations.markings,arrows}

\usepackage{wrapfig}
\usepackage{color}
\usepackage{graphicx,framed}
\usepackage[colorlinks=true, pdfstartview=FitV, linkcolor=darkblue, citecolor=darkblue, urlcolor=darkblue]{hyperref}
\definecolor{shadecolor}{rgb}{0.95, 0.95, 0.86}
\definecolor{darkgreen}{rgb}{0.2, 0.5,  0}

\def\&{\vspace{-5pt}&}

\def \C{\mathbb C}

\def \eqref#1{(\ref{#1})}
\def \& {&\hspace{-10pt}}

\newcommand{\bt}{\beta}

\renewcommand{\O}{\Omega}

\newcommand{\br}{{\mathbb R}}

     \newcommand{\lop}[1]{\mathfrak{L}(#1)}
      \newcommand{\Lie}{\mathrm{Lie}}

\newtheorem{theorem}{Theorem}[section]
\newtheorem{example}[theorem]{Example}
\newtheorem{exercise}[theorem]{Exercise}

\newtheorem{lemma}[theorem]{Lemma}
\newtheorem{remark}[theorem]{Remark}

\newtheorem{proposition}[theorem]{Proposition}
\newtheorem{corollary}[theorem]{Corollary}
\newtheorem{definition}[theorem]{Definition}

\def\bt{\begin{theorem}}
\def\et{\end{theorem}}
\def\bc{\begin{corollary}}
\def\ec{\end{corollary}}
\def\bx{\begin{example}}
\def\ex{\end{example}}
\def\bxr{\begin{exercise}\small}
\def\exr{\end{exercise}}
\def\bl{\begin{lemma}}
\def\el{\end{lemma}}
\def\bd{\begin{definition}}
\def\ed{\end{definition}}
\def\bp{\begin{proposition}}
\def\ep{\end{proposition}}

\def\br{\begin{remark}}
\def\er{\end{remark}}

\def\be{\begin{equation}}
\def\ee{\end{equation}}

\def\&{\hspace{-15pt}&}
\def\bea{\begin{eqnarray}}
\def\eea{\end{eqnarray}}

\def\L{\mathrm {Lie}}

\def\1{{\bf 1}}

\def\W{{\mathcal W}}

\newcommand{\F}{\mathbb F}

\makeatletter
\@addtoreset{equation}{section}
\makeatother

\begin{document}
\title{Conjugate Frobenius manifold  and inversion symmetry}

\maketitle
\begin{center}
Zainab Al-Maamari
\footnote{Sultan Qaboos University, Muscat, Oman, s100108@student.squ.edu.om},
Yassir Dinar
\footnote{Sultan Qaboos University, Muscat, Oman, dinar@squ.edu.om}.
\end{center}

\maketitle

\begin{abstract}
We give a conjugacy relation  on  certain type of Frobenius manifold structures using the theory of flat pencils of metrics. It leads to a geometric interpretation for the inversion symmetry of solutions to  Witten-Dijkgraaf-Verlinde-Verlinde (WDVV) equations. 
\end{abstract}

\maketitle
\tableofcontents

\section{Introduction}

Boris Dubrovin introduced the notion of a  Frobenius manifold as a geometric realization of a potential $\mathbb F$  {which} satisfies a system of partial differential equations known in topological field theory as Witten-Dijkgraaf-Verlinde-Verlinde (WDVV) equations. More precisely, a Frobenius algebra is a commutative associative algebra with an identity $e$ and a nondegenerate bilinear form $\Pi$  {compatible with the product}, i.e., $\Pi(a\circ b,c)=\Pi(a,b\circ c)$.
A Frobenius manifold is a manifold  with
a smooth structure of a Frobenius algebra on the tangent space  at any point  with certain compatibility conditions. Globally, we  require the metric $\Pi$ to be flat and the identity vector field $e$ to be covariantly constant with respect to the corresponding Levi-Civita connection. Detailed information about Frobenius manifolds and related topics can be found in $\cite{DuRev}$.

Let $M$ be a Frobenius manifold. In flat  coordinates $(t^1,...,t^r)$ for $\Pi$ where $e= \partial_{t^{r}}$ the compatibility conditions imply  that there exists a function $\mathbb{F}(t^1,...,t^r)$ which encodes  the Frobenius structure, i.e., the flat metric is given by
\be \label{flat metric} \Pi_{ij}(t)=\Pi(\partial_{t^i},\partial_{t^j})=  \partial_{t^{r}}
\partial_{t^i}
\partial_{t^j} \mathbb{F}(t)\ee
and, setting $\Omega_1(t)$ to be the inverse of the matrix $\Pi(t)$, the structure constants of the Frobenius algebra are given by 
\[ C_{ij}^k(t)=\Omega_1^{kp}(t)  \partial_{t^p}\partial_{t^i}\partial_{t^j} \mathbb{F}(t).\]
 Here, and in what follows, summation with respect to repeated
upper and lower indices is assumed. The  definition includes the existence  of a vector field $E$ of the form $E=(a_i^j t^i+b^j)\partial_{t^j}$ satisfying
\be  \label{quasihomog1}
E\mathbb F(t)= \left(3-d \right) \mathbb{F}(t)+ \frac{1}{2}A_{ij} t^i t^j+B_i t^i+c
\ee 
where  $a_i^j$, $b_j$, $c$, $A_{ij}$,  $B_i$ and $d$ are constants with $a_r^r=1$. The vector field $E$ is called the Euler vector field and the number $d$ is called the charge of the Frobenius manifold. 
The  associativity  of  Frobenius
algebra implies that the potential $\mathbb{F}(t)$ satisfies   the WDVV  equations
\begin{equation} \label{frob}
 \partial_{t^i}
\partial_{t^j}
\partial_{t^k} \mathbb{F}(t)~ \Omega_1^{kp} ~\partial_{t^p}
\partial_{t^q}
\partial_{t^n} \mathbb{F}(t) = \partial_{t^n}
\partial_{t^j}
\partial_{t^k} \mathbb{F}(t) ~\Omega_1^{kp}~\partial_{t^p}
\partial_{t^q}
\partial_{t^i} \mathbb{F}(t),~~ \forall i,j,q,n. 
  \end{equation}
Conversely, an arbitrary  potential $\mathbb F(t^1,\ldots, t^r)$ satisfying equations \eqref{frob} and \eqref{quasihomog1} with \eqref{flat metric} determines a Frobenius manifold structure  on its domain \cite{DuRev}. Moreover, there exists a quasihomogenius flat pencil of metrics (QFPM) of degree  $d$ associated to  the Frobenius structure on $M$  which consists of the intersection form ${\Omega}_2$ and the flat metric $\O_1$ with the function $\tau=\Pi_{i1}t^i$ (see definition \ref{FPM} below). Here 
\be \label{intersection form}
\O_2^{ij}(t):=E(dt^i\circ dt^j)
\ee
where the product $dt^i\circ dt^j$ is defined by lifting the product on $TM$ to $T^*M$ using the flat metric $\O_1$.
  In this article we prove that, when  $d\neq 1$, $e(\tau)=0$ and $E(\tau)=(1-d)\tau$, we can construct another  QFPM of degree $2-d$ on $M$ consisting  of the intersection form $\O_2$ and a different flat metric $\widetilde\O_1$. We call it the conjugate QFPM.  In particular, under a specific regularity condition, we get a conjugation  between a certain type of Frobenius manifold structures on a given manifold. Precisely, we  prove the following theorem. 

\bt \label{dual Frob manif}
Let $M$ be a Frobenius manifold with the Euler vector field $E$ and the identity vector field $e$. Suppose the associated QFPM  is regular of degree $d$ with  a function $\tau$. Assume that 
$e(\tau)=0$ and $E(\tau)=(1-d)\tau$. Then  we can construct another  Frobenius manifold structure on $M\backslash \{\tau=0\}$ of degree $2-d$. Moreover, we can apply the same  method to the new Frobenius manifold structure  and it leads to the original Frobenius manifold structure. 
\et  
For a fixed Frobenius manifold  the new structure that can be obtained using Theorem \ref{dual Frob manif} will be called the conjugate Frobenius manifold structure.

 Let us assume  $\Pi_{i,j}=\delta_{i+j}^{r+1}$, i.e., the potential  $\F$ has the standard form 
\be \label{norm potential}
 \F(t) = \frac{1}{2} (t^r)^2 t^1 + \frac{1}{2} t^r \sum_{i=2}^{r-1} t^i t^{r-i+1} + G(t^1,...,t^{r-1})
\ee 
and the quasihomogeneity condition \eqref{quasihomog1} takes the form
\begin{equation}\label{quasihomog}  E= d_i t^i\partial_{t^i},~ E\mathbb F(t)= \left(3-d \right) \mathbb{F}(t);~~d_{r}=1.\end{equation} 
Here, the numbers $d_i$ are called the degrees of the Frobenius manifold. Recall that a symmetry of  the WDVV equations is a transformation of the form  
\[ 
t^i\mapsto z^i, ~ {\Pi}\mapsto \widetilde{\Pi},~ \mathbb F\mapsto \widetilde  \F
\]
such that $\widetilde \F$ satisfies  the WDVV equations. 
The inversion symmetry (\cite{DuRev}, Appendix B)  is an involutive  symmetry given by setting 
\be\label{Dob coord}
z^1=-\frac{1}{t^1},~z^r=\Pi_{ij}(t)\frac{t^i t^j}{2 t^1},~ z^k=\frac{t^k}{t^1}, ~2\leq k< r.
\ee 
Then 
\be\label{inv potential}
\widetilde \F(z) :=(t^1)^{-2}\left( \F(t)-\frac{1}{2}t^r \Pi_{ij}t^i t^j\right) \ee 
is another solution to the WDVV equations with the flat metric $\widetilde\Pi_{ij}(z)=\delta_{i+j}^{r+1}$. The   charge of the corresponding Frobenius manifold structure is $2-d$ and the degrees are 
\be \label{degrees inv} \widetilde{d}_1=-d_1,\ \ \widetilde{d}_r=1, \ \ \widetilde{d}_i = d_i-d_1 \ \ for \ \ 1< i < r.\ee
The inversion  symmetry is obtained from a special Schlesinger
transformation of the system of linear ODEs with rational
coefficients associated to the WDVV equations. A   geometric  relation between Frobenius manifold structures correspond to $\F(t)$ and $\widetilde \F (z)$ was outlined  through the sophisticated notion of Givental groups in \cite{givental}. In this article, we obtained a simple  geometric  interpertation and we report that $\widetilde \F(z)$ is the potential of the conjugate Frobenius manifold structure. In other words, we prove the following theorem.  

\bt \label{main thm}
Let $M$ be a Frobenius manifold with charge $d\neq 1$. Suppose in the flat coordinates $(t^1,\ldots,t^r)$, the potential $\mathbb F(t)$ has the standard form \eqref{norm potential} and the quasihomogeneity condition takes the form  \eqref{quasihomog} with $d_i\neq \dfrac{d_1}{2}$ for every $i$.  Then we can construct the  conjugate   Frobenius manifold structure on $M\backslash\{t^1=0\}$.    Moreover,  flat coordinates  for the conjugate Frobenius manifold  are 
\be \label{change coord}
s^1= -t^1 , \ \ s^i= t^i (t^1)^{\frac{d_1-2d_i}{d_1}}\ \ for \ \ 1<i< r, \ \ s^r= \frac{1}{2} \sum_{i=1}^{r}  t^i t^{r-i+1} (t^1)^{\frac{-2}{d_1}-1}.
\ee 
In addition, the corresponding potential   equals the potential  obtained by applying the inversion symmetry to $\mathbb F(t)$  and it is given by   \be
\widetilde{\F}(s) = (t^1)^{ \frac{-4}{d_1}}  \left(\F(t^1,\ldots,t^r) -\frac{1}{2} t^r \sum_{1}^{r} t^i t^{r-i+1}\right). 
\ee 
\et

Examples of Frobenius manifolds satisfying the hypotheses of Theorem \ref{main thm} include Frobenius manifold structures constructed on orbits spaces of standard reflection representations of irreducible Coxeter groups in \cite{DCG} and \cite{polyZuo} and algebraic Frobenius manifolds constructed using classical $W$-algebras \cite{mypaper5}. However, the result presented in this article is a consequence of the  work \cite{dicy} and \cite{nonref}. There, we investigated the existence of Frobenius manifold structures on orbits spaces of some non-reflection representations of  finite groups and we noticed that certain structures appear in pairs.  Analyzing such pairs led us to the notion of   conjugate Frobenius manifold.  
 

This article is organized as follows. In section \ref{flat frob}, we review the relation between Frobenius manifold, flat pencil of metrics and  compatible Poisson brackets of hydrodynamic type. Then we introduce a conjugacy  relation between  certain class of  quasihomogeneous flat pencils of metrics in section \ref{dualityFrob}. It can be interpreted  as a conjugacy relation between certain class of compatible Poisson brackets of hydrodynamic type. We prove Theorem \ref{dual Frob manif} in section \ref{dualityFrob} and Theorem \ref{main thm} in  section \ref{relation F and new F}. In  section \ref{first section}, we discuss the findings of this article on polynomial Frobenius manifolds. We end the article with some  remarks.


\section{Background} \label{flat frob}

We  review in this section  the relation between flat pencil of metrics,  compatible Poisson brackets of hydrodynamics type and  Frobenius manifold. More details can be found in \cite{Du98}.

Let $M$ be a smooth manifold of dimension $r$ and  fix  local coordinates $(u^1, . . . , u^r)$ on $M$. 

\bd \label{contra metric} A symmetric bilinear form $(. ,. )$ on $T^*M$ is called a contravariant
metric if it is invertible on an open dense subset $M_0 \subseteq M$. We define the contravariant Christoffel symbols $\Gamma^{ij}_k$  for a contravariant
metric $(. ,. )$ by
\[
\Gamma^{ij}_k:=-\Omega^{im} \Gamma_{mk}^j
\]
where $\Gamma_{mk}^j$ are the  Christoffel symbols of the metric $<. ,. >$ defined on $TM_0$ by the inverse of the matrix $\Omega^{ij}(u)=(du^i, du^j)$.
We say the metric $(.,.)$ is flat if  $<. ,. >$ is flat.
\ed

Let $(. ,. )$ be a contraviariant metric on $M$ and set $\O^{ij}(u)=(du^i, du^j)$. Then we will use $\Omega$ to refer to the metric and $\Omega(u)$ to refer to its  matrix  in the coordinates. In particular, the Lie derivative of $(. ,. )$ along a vector field $X$  will be written $\Lie_X \Omega$ while $X\Omega^{ij}$ means the vector field $X$ acting on the entry $\Omega^{ij}$. The Christoffel symbols given in definition \ref{contra metric} determine for $\O$ the contravariant (resp. covariant) derivative  $\nabla^{i}$ (resp.  $\nabla_{i}$) along the covector $du^i$ (resp. the vector field $\partial_{u^i}$). They are related by the identity  $\nabla^{i}=\O^{ij}(u) \nabla_{j}$.

\bd 
A flat pencil of metrics (FPM)  on $M$ is a pair $(\Omega_2,\Omega_1)$ of 
 two flat contravariant metrics $\O_2$ and $\O_1$ on $M$ satisfying 
 \begin{enumerate}
     \item $\O_2+\lambda \O_1$ defines a flat metric on $T^*M$ for  a generic constant $\lambda$,
     \item the Christoffel symbols of $\O_2+\lambda \O_1$ are   $\Gamma_{2k}^{ij}+\lambda \Gamma_{1k}^{ij}$, where   $\Gamma_{2k}^{ij}$ and $ \Gamma_{1k}^{ij}$ are the Christoffel symbols of $\O_2$ and $\O_1$,  respectively. 
 \end{enumerate}  
\ed  
\bd \label{FPM}  A  flat pencil of metrics $(\O_2,\O_1)$ on  $M$ is called quasihomogeneous flat pencil of metrics (QFPM) of  degree $d$ if there exists a function $\tau$ on $M$ such that the 
 vector fields $E$ and $e$ defined by 
\begin{eqnarray} \label{tau flat pencil} E&=& \nabla_2 \tau, ~~E^i
=\O_2^{ij}(u)\partial_{u^j}\tau
\\\nonumber  e&=&\nabla_1 \tau, ~~e^i
= \O_1^{ij}(u)\partial_{u^j}\tau  \end{eqnarray} satisfy 
\be \label{vector fields} [e,E]=e,~~ \Lie_E \O_2 =(d-1) \O_2,~~ \Lie_e \O_2 =
\O_1~~\mathrm{and}~~ \Lie_e\O_1
=0.
\ee 
Such a QFPM is  \textbf{regular} if  the
(1,1)-tensor
\begin{equation}\label{regcond}
  R_i^j = \frac{d-1}{2}\delta_i^j + {\nabla_1}_i
E^j
\end{equation}
is  nondegenerate on $M$.
\ed 

Let $(\O_2,\O_1)$ be a QFPM of degree $d$. Then according to \cite{Du98}, we can  fix flat coordinates $(t^1,t^2,\ldots,t^r)$ for $\O_1$ such that 
\be \label{gamma id}
\tau=t^1, \ E^i= \Omega_2^{i1}, \ e^i= \Omega_1^{i1}, \ \ \Gamma_{1,k}^{ij}=0, \ \
\Gamma_{2,k}^{i1}= \frac{1-d}{2} \delta_k^i, \ \
\Gamma_{2,k}^{1j}= \frac{d-1}{2}\delta_k^j+ \partial_{t^k} E^j,\ \ \partial_{t^1} E^1=1-d.
\ee
Moreover, if $(\O_2,\O_1)$ is regular then $d\neq 1$.

 Consider the loop space $\lop M$ of $M$, i.e., the space of smooth maps from the circle $S^1$ to $M$. A local Poisson bracket  on $\lop M$ is a Lie algebra structure on the space of local functionals on $\lop M$. Let $\{.,.\}$ be a local Poisson  bracket of hydrodynamic type (PBHT), i.e., it has the following form in the local coordinates  \cite{Du98}
\be\label{poisson} \{u^i(x),u^j(y)\}= \Omega^{ij}(u(x))\delta' (x - y) + \Gamma_k^{ij} (u(x)) u_x^k  \delta (x-y), \, i,j=1,\ldots,r\ee 
where  $\delta(x-y)$ is the Dirac  delta function defined by
$\int_{S^1} f(y) \delta(x-y) dy=f(x)$. Then we say $\{.,.\}$  is nondegenerate   if  $\det \Omega^{ij}\neq 0$ and the Lie derivative of $\{.,.\}$ along  a vector field $X:=X^i\partial_{u^i}$  reads 
\begin{align*}
\Lie_X\{.,.\}(u^i(x),u^j(y))&= (X^s\partial_{u^s}\Omega^{ij}- \Omega^{s j}\partial_{u^s}X^{i}-\Omega^{is}\partial_{u^s} X^{j})\delta'(x-y)
\\\nonumber &+(X^s \partial_{u^s}\Gamma_k^{ij}-\Gamma_k^{sj} \partial_{u^s}X^i-\Gamma^{i s}_k\partial_{u^s}X^j+\Gamma_s^{i j} \partial_{u^k}X^s-\Omega^{i s}\partial_{u^s}\partial_{u^k} X^j)u_x^k\delta(x–y).
\end{align*}
We will use the following two theorems. 
\bt\label{serg}\cite{serg}
Let $X$ be a vector field on $M$ and $\{ . , .\}$ be a PBHT on $\lop M$. If $\Lie_X^2 \{.,.\}=0$, then $\Lie_X \{.,.\}$ is a PBHT and it is compatible with $\{.,.\}$, i.e., $\{.,.\} + \lambda \Lie_X \{. , .\}$ is a PBHT for every constant $\lambda$. 
\et

\bt\label{Du and Nov}\cite{Du and Nov}
The form \eqref{poisson} defines a  nondegenerate PBHT $\{ . , .\}$ if and only if the matrix  $\Omega^{i j}(u)$ defines a flat contravariant metric on $M$ and $\Gamma_k^{i j}(u)$ are its Christoffel symbols.
\et

 From Theorem \ref{Du and Nov} and Theorem \ref{serg}, we get the following corollary: 

\bc \label{FPM and PBHT}
Let $\{.,.\}_2$ and $\{.,.\}_1$ be two nondegenerate compatible PBHT on $\lop M$ having the form 
\[ \{u^i(x),u^j(y)\}_\alpha = \Omega_\alpha^{ij} (u(x))\delta' (x-y) + \Gamma_{\alpha,k}^{ij} (u(x)) u_x^k  \delta (x – y),~ \alpha=1,2.\]
Suppose $\{.,.\}_2+\lambda \{.,.\}_1$ is a nondegenerate PBHT for a generic constant $\lambda$. Then $(\Omega_2,\Omega_1)$ is a FPM on $M$. Conversely, a FPM on $M$ determines  nondegenerate compatible Poisson brackets of hydrodynamic type on $\lop M$.
\ec

As mentioned in the introduction, if $M$ is a Frobenius manifold of charge $d$ then there is an associated QFPM $(\O_2,\O_1)$ of degree $d$ on $M$, where $\O_2$ is the intersection form and $\O_1$ is the flat metric. In the flat coordinates $(t^1,\ldots,t^r)$ we have   $\tau= \Pi_{i 1} t^i$. Then   the Euler vector field $E$ and the identity vector field $e$ of the Frobenius manifold   have the form  \eqref{tau flat pencil}  and satisfy equations \eqref{vector fields}. The following theorem give a converse statement. 
\bt\cite{Du98}\label{dub flat pencil} Let $M$ be a manifold carrying a regular QFPM  $(\Omega_2,\Omega_1)$ of degree $d$. Then there exists a unique Frobenius manifold structure on $M$ of charge $d$ where $(\Omega_2,\Omega_1)$ is the associated QFPM. 
\et

\section{Conjugate Frobenius manifold} \label{dualityFrob}

We fix a   manifold $M$ with a QFPM $T=(\Omega_2,\Omega_1)$ of degree $d\neq 1$. We fix a function $\tau$ for $T$ which determines  the vector fields $E$ and $e$  (see definition \ref{FPM}). We suppose 
\be\label{new condition} e(\tau)=0 \ \ \textrm{and} \ \ E(\tau)= (1-d) \tau.\ee
We introduce the function  $f(\tau):=(\tau)^{\frac{2}{1-d}}$ and  the vector field $ \widetilde{e} := f(\tau) e$.  We define \be 
\widetilde{\Omega}_1: = \L_{\widetilde{e}} \Omega_2 = f \Omega_1 -f'(E \otimes e + e \otimes E).
\ee 
Then 
\begin{align}\label{new o}
    \L^2_{ \widetilde{e}} {\Omega}_2 & =  f^2 (\L_e^2 \Omega_2) + (2(f')^2 E(\tau) -4f f')e\otimes e +f f'e(\tau) \Omega_1\\ \nonumber  &+  ((f')^2- f f'') e(\tau)(E\otimes e+e\otimes E)=0
\end{align}

We fix  flat coordinates $(t^1,\ldots,t^r)$ leading to the identities \eqref{gamma id}. Considering the condition \eqref{new condition}, we will further assume  that $e=\partial_{t^r}$. Thus
\be
 \Omega_1^{i1}=\delta^{i}_r, \ \ \partial_{t^r} \Omega_2^{i1}=\partial_{t^r}E^i=\delta^{i}_r.
\ee 

 Let $\{.,.\}$ denote the nondegenerate PBHT associated to $\Omega_2$. Then by  Corollary \ref{FPM and PBHT}, $\L_{e}\{.,.\}$ is the PBHT associated to $\O_1$ and  $\L_{e}^2\{.,.\}=0$. We have a similar statement for $\widetilde e$. 
\bp\label{PBHT}
$\L_{\widetilde{e}}^2 \{.,.\}=0$. In particular,  $\L_{\widetilde{e}} \{.,.\}$ is a PBHT compatible with $\{.,.\}$.
\ep
\begin{proof}
  The PBHT associated to $\Omega_2$ has the form 
  \[ \{t^\alpha(x),t^\beta(y)\} =\Omega_2^{\alpha \beta} \delta' (x-y) + \Gamma_{2,\gamma}^{\alpha \beta}  t^\gamma_x \delta(x-y).\]
Here and in what follows, it is to be understood that all functions on the right hand side depend on $t(x)$.   Note that 
  \[ 
  \L_{\widetilde e}\{.,.\}(t^\alpha(x),t^\beta(y))=\widetilde \O_1^{\alpha\beta}\delta'(x-y)+\widetilde\Gamma^{\alpha \beta}_{2,\gamma} t_x^\gamma \delta(x-y)
\]
where 
\begin{align*}
\widetilde \Gamma_{2,\gamma}^{\alpha \beta}& = { \widetilde{e}}^\varepsilon \partial_\varepsilon \Gamma_{2,\gamma}^{\alpha \beta}- \Gamma_{2,\gamma}^{\varepsilon \beta} \partial_\varepsilon { \widetilde{e}}^\alpha – \Gamma_{2,\gamma}^{\alpha \varepsilon} \partial_\varepsilon { \widetilde{e}}^\beta + \Gamma_{2,\varepsilon}^{\alpha \beta} \partial_\gamma { \widetilde{e}}^\varepsilon – \Omega^{\alpha \varepsilon}_2 \partial_{\varepsilon \gamma}^2 { \widetilde{e}}^\beta\\
 &   = - \Gamma_{2,\gamma}^{\varepsilon \beta} \delta_r^\alpha \delta_\varepsilon^1 f' – \Gamma_{2,\gamma}^{\alpha \varepsilon} \delta_r^\beta \delta_\varepsilon^1 f' + \Gamma_{2,\varepsilon}^{\alpha \beta} \delta_r^\varepsilon \delta_\gamma^1 f' – \Omega^{\alpha \varepsilon}_2 \delta_r^\beta \delta_\gamma^1 \delta_\varepsilon^1 f''.
\end{align*}
From equation \eqref{new o}, the coefficients of $\delta'(x-y)$ of  $\L_{ \widetilde{e}}^2 \{.,.\}$ vanish while   the coefficients  $\widetilde{\widetilde{\Gamma}}_{2,\gamma}^{\alpha \beta}$ of $\delta(x-y)$ have the form 
\begin{align*}
 \widetilde{\widetilde{\Gamma}}_{2,\gamma}^{\alpha \beta}=& -f f'' \partial_r \Omega_2^{\alpha \varepsilon} \delta_r^\beta \delta_\gamma^1 \delta_\varepsilon^1 + f'^2 \delta_r^\alpha \delta_r^\beta \delta_m^1 \delta_\varepsilon^1 \Gamma_{2,\gamma}^{m \varepsilon} - f'^2 \delta_r^\beta \delta_r^m \delta_\gamma^1 \delta_\varepsilon^1 \Gamma_{2,m}^{\alpha \varepsilon } \\
& + f'^2 \delta_r^\beta \delta_r^\alpha \delta_\varepsilon^1 \delta_m^1 \Gamma_{2,\gamma}^{\varepsilon m }- f'^2 \delta_r^\alpha \delta_r^m \delta_\gamma^1  \delta_\varepsilon^1 \Gamma_{2,m}^{\varepsilon \beta}+ f' f'' \Omega_2^{\varepsilon m}\delta_\varepsilon^1 \delta_r^\alpha \delta_r^\beta  \delta_\gamma^1 \delta_m^1\\
& - f'^2 \delta_\gamma^1 \delta_r^\beta \delta_m^1 \delta^\varepsilon_r \Gamma_{2,\varepsilon}^{\alpha m}
- f'^2 \delta_\gamma^1 \delta_r^\alpha \delta_m^1 \delta^\varepsilon_r \Gamma_{2,\varepsilon}^{ m \beta} -\widetilde{\Omega}_2^{\alpha \varepsilon} \delta_r^\beta \delta_\gamma^1  \delta_\varepsilon^1 f''.
\end{align*}
Then from the identities \eqref{gamma id} and the definition of $f(\tau)$, it follows that  $\widetilde{\widetilde{\Gamma}}_{2,\gamma}^{\alpha \beta}=0$. For example,
 \begin{align*}
\widetilde{\widetilde{\Gamma}}_{2,1}^{r r}& = - f \partial_r{\Omega_2^{r 1} f''} + f'^2 \Gamma_{2,1}^{1 1}- f'^2 \Gamma_{2,r}^{1 r}+ \Omega_2^{1 1} f'' f'+ f'^2 \Gamma_{2,1}^{1 1}- f'^2 \Gamma_{2,r}^{r 1}- f'^2 \Gamma_{2,r}^{r 1} - f'^2 \Gamma_{2,r}^{r 1}- \widetilde{\Omega}_1^{r 1} f''\\
&= -(d+1) f'^2 +(1-d) \tau f'f'' -f f'' -(-f) f'=0
\end{align*}  and when  $\gamma =1$, $\alpha=r$ and $\beta \neq r$ 
\begin{align*}
\widetilde{\widetilde{\Gamma}}_{2,1}^{r \beta}&=  - 2 f'^2   \Gamma_{2,r}^{1 \beta}= - 2 f'^2 (\frac{d-1}{2}\delta^\beta_r+ \partial_{t^r} E^\beta)=0.\end{align*}\end{proof}

\bl\label{new flat metrics} The pair $\widetilde T=(\O_2,\widetilde{\Omega}_1)$  form a QFPM of degree $\widetilde{d}=2-d$. Moreover, if $T$ is regular then $\widetilde T$ is regular.
\el 
\begin{proof} The second term of the identity  \[\widetilde{\Omega}_1(t)= f \Omega_1 -f' E^i (\partial_{t^i} \otimes \partial_{t^r} + \partial_{t^r}  \otimes  \partial_{t^i})\] contributes only  to entries of the last row and last column of $ \widetilde{\Omega}_1(t)$.
From the normalization of $\Omega_1$, we get \[\widetilde{ \Omega}_1^{i1}(t)=(f-f' E(\tau))
 \delta^i_r= (f- (1-d) \tau  f') \delta^i_r= (- f) \delta^i_r. \]
Therefore,
\[    \det  \widetilde{\Omega}_1(t)= f^r \det \Omega_1(t) \neq 0.
\]
Hence, using Proposition \ref{PBHT} and Corollary \ref{FPM and PBHT}, $\widetilde T$ is a FPM. Let $\widetilde \nabla$ denote the contravariant (and also the covariant) derivative of $\widetilde \O_1$ and  set $ \widetilde{\tau}:=-\tau=-t^1$. Then the vector fields  \[\widetilde{e}:=\widetilde{\nabla}_1   \widetilde{\tau}, ~\text{and} \ \ \widetilde{E}:={\nabla}_2 \widetilde{\tau}=-E\]
 satisfy equations \eqref{vector fields}  and   \be
\L_{\widetilde{E}} \Omega_2 = \L_{-E} \Omega_2= -(d-1)\Omega_2= (\widetilde{d}-1) \Omega_2.\ee
Hence,  $\widetilde T$  is a QFPM  of degree $\widetilde{d}=2-d$. For the regularity condition \eqref{regcond}, we have
 \be
  \widetilde R_i^j(t) = \frac{\widetilde{d}-1}{2}\delta_i^j + \widetilde{\nabla}_{1i}
(-E^j)=\frac{1-d}{ 2} \delta_i^j-{\nabla}_{1 i}
(E^j) =- R_i^j(t).
\ee
 Therefore, $\det (\widetilde R_i^j) \neq 0$ if and only if $\det (R_i^j) \neq 0$. 
\end{proof}

We keep the definitions    $\widetilde \tau=-\tau$ and  $\widetilde E=-E$  given in the proof of Lemma \ref{new flat metrics} and we call  $\widetilde T=(\O_2,\widetilde{\Omega}_1)$ the conjugate QFPM of $T$. The name is motivated by the following corollary. 
\bc\label{cor on duality}
$\widetilde T$ has a conjugate and it equals  $T$.
\ec 

\begin{proof}
We observe that $\widetilde d=2-d\neq 1$ and the function  $\widetilde{\tau}=-\tau$ satisfies the requirements \eqref{new condition} as  
\be   \widetilde{e} (\widetilde{\tau})=0 \ \ \text{and} \ \  \widetilde{E} ( \widetilde{\tau})= -E (-t^1)=(1-d) t^1= (1- \widetilde{d}) \widetilde{\tau}.\ee 
However, applying Lemma \ref{new flat metrics} to $\widetilde T$, we get a QFPM $(\O_2,\L_{\widetilde{\widetilde{e}}}\O_2)$ where 
\[     \widetilde{ \widetilde{e}}=f( \widetilde{\tau})  \widetilde{e}= \widetilde{\tau}^{\frac{2}{1- \widetilde{d}}}  \,\widetilde{e}= (t^1)^{\frac{2}{1- \widetilde{d}}}.(t^1)^{\frac{2}{1-{d}}} \partial_{t^r}=e.
\]
\end{proof}

Now we can prove Theorem \ref{dual Frob manif}.

\begin{proof}[Proof of Theorem  \ref{dual Frob manif}]
From the work in \cite{Du98},  regularity of the associated QFPM implies that  the charge  $d\neq 1$. Then the proof follows from applying Lemma \ref{new flat metrics}, Corollary \ref{cor on duality} and Theorem \ref{dub flat pencil} to the associated regular QFPM.
\end{proof}
For a fixed  Frobenius manifold, the new Frobenius manifold structure  constructed using Theorem \ref{dual Frob manif} will be called the conjugate Frobenius manifold structure. 
\bx \label{example dim 2}
We consider the Frobenius manifold structure of charge -1 defined by the following solution  to the WDVV equations.
\[ 
\F=\frac{1}{2} t_2^2 t_1 + t_1^2\log t_1
\]
In the examples, we use subscript indices instead of superscript indices for convenience. Here, the identity vector field  $e=\partial_{t_2}$ and  the Euler vector field   $E=2 t_1 \partial_{t_1}+t_2 \partial_{t_2}$. Note that    $EF=(3-d) F + 2 t_1^2$. 
The corresponding regular QFPM consists of  
\be 
\Omega_2(t) =\left(
\begin{array}{cc}
 2 t_1 & t_2  \\
 t_2 & 4 \\
\end{array}
\right),~\Omega_1(t) =\left(
\begin{array}{cc}
 0 & 1 \\
 1 & 0 \\
\end{array}
\right).
\ee 
The conjugate QFPM $\widetilde T=(\O_2,\widetilde \O_1)$ is of degree $\widetilde d=3$. In the  coordinates 
\[
s_1=-t_1,\ \ s_2=\frac{t_2}{t_1} \] 
we have 
\[\O_2(s)=\left(
\begin{array}{ccc}
 -2s_1 & s_2  \\
 s_2 & \frac{4}{s_1^2}\\
\end{array}
\right), ~\widetilde\O_1(s)=\left(
\begin{array}{ccc}
 0 & 1  \\
1 & 0\\
\end{array}
\right)\]
and the potential of the conjugate Frobenius manifold structure has the form \[ \widetilde{\F}=\frac{1}{2} s_1 s_2^2- \log s_1. \]
Note that the Euler vector field  $\widetilde{E}=-E(s)=-2s_1 \partial_{s_1}+s_2 \partial_{s_2}$ and $\widetilde E \widetilde\F=(3-\widetilde d)\widetilde\F+2$. We observe that applying the inversion symmetry to the potential $\F(t)$, we get 
\[\widehat\F(z)=\frac{1}{2} z_1 z_2^2- \log z_1+~ \text{constant} \]
and  $\widehat\F(z)$  defines  the same  conjugate Frobenius manifold structure. We prove   this for certain type of Frobenius manifolds  in next section. 
\ex

\bx\label{non regular but has new F}
We consider Frobenius manifold structures found recently in  \cite{new ref} on the orbits space of the reflection group of type $B_4$. It is provided to us by the anonymous reviewer of this article as an example of Frobenius manifold structure  whose  associated QFPM has a conjugate but it is not regular.  The  potential of this Frobenius manifold reads
$$\F=\frac{1}{2} t_4^2 t_1+t_2 t_3 t_4-\frac{1}{72}t_1^4+\frac{1}{2} t_3 t_1^2+\frac{1}{6} t_2^2 t_3 t_1-\frac{9 }{4}t_3^2+\frac{1}{108} t_2^4 t_3+\frac{3}{2}
   t_3^2 \log t_3.$$
where the charge and degrees given by
$$d=\frac{1}{3}, \ \ d_1=\frac{2}{3},\ \  d_2=\frac{1}{3},\ \  d_3=\frac{4}{3},\ \ d_4=1.$$
The action of the Euler vector field reads 
\be \label{quasi not} E\ \F (t)= \left(3-d \right) \F (t)+ \frac{1}{2}A_{ij} t^i t^j=\left(3-d \right) \F (t)+2 t_3^2 \ee 
and  the intersection metric $\O_2$ will be
\be\label{get O dif}
{\dot{\Omega}}_2^{i j}(t)= \O_2^{ij}(t)+ A^{i j},~~A^{i j}= {\O}_1^{i\alpha}(t) {\O}_1^{j \beta}(t) A_{\alpha \beta}.
\ee
The associated QFPM $T=(\O_2,\O_1)$ is  not regular. However,  it has a conjugate QFPM $\widetilde T=(\O_2,\widetilde \O_1)$. Flat coordinates $(s_1,s_2,s_3,s_4)$ for $\widetilde \O_1$ are defined by 
\[
t_1= -s_1,\ \ t_2= s_2,\ \ t_3= -s_1^3 s_3,\ \ t_4= -s_4 s_1^3-s_2 s_3 s_1^2\] 
Note that one can still apply the inversion symmetry to the potential $\mathbb F(t)$ to get a Frobenius manifold structure with a potential $\widehat  \F(z)$ \cite{DuRev}. We checked that  the  QFPM obtained from $\widehat  \F(z)$ agrees with $\widetilde T$. We do not consider this type of Frobenius manifolds in the next section as we will assume regularity condition \eqref{regcond} of the quasihomogeneous  flat pencils of metrics. 
\ex

Let us assume $E$ has the form  $E=d_i t^i \partial_{t^i}$.  Then $d_1=1-d$ and we have the following standard results.

\bc \label{reg bcz E}
$T$ is regular QFPM if and only if $d_i \neq \frac{d_1}{2}$ for all $i$.
\ec
\begin{proof} Applying the definition \ref{FPM} to the  matrix  $R_i^j(t)=(\frac{d-1}{2}+d_i) \delta_i^j=(-\frac{d_1}{ 2}+d_i) \delta_i^j$.
\end{proof}

\bl\label{degrees}
If $\Omega_1^{i j} \neq 0$, then $d_i+d_j =2-d$. Thus, if the numbers $d_i$ are all distinct then we can choose the coordinates $(t^1,\ldots,t^r)$ such that $\Omega_1^{i j}=\delta^{i+j}_{r+1}$.
\el
\begin{proof} Notice that using $[e,E]=e$, we get $\L_E \Omega_1 = (d-2) \Omega_1$. Then the statement follows from  the equation 
\[ (2-d) \Omega_1^{ij}(t)=\L_E \Omega_1^{ij} (dt^i,dt^j) = -d_i \Omega_1 (dt^i,dt^j)-d_j \Omega_1 (dt^i,dt^j).\]
 \end{proof}

\section{Relation with inversion symmetry} \label{relation F and new F}
We continue using  notations and assumptions  given in the previous section, but we suppose that $T$ is regular. Consider the Frobenius manifold structure defined on $M$ by Theorem \ref{dub flat pencil} and let $\F(t)$ be the corresponding potential. We assume  $\Omega^{ij}_1(t)=\delta^{i+j}_{r+1}$ which is  equivalent to requiring  that  $\F(t)$  has the standard form \eqref{norm potential}. We suppose further that the quasihomogeneity condition for $\F(t)$ takes the form \eqref{quasihomog}. In this case the intersection form $\O_2$ satisfies \cite{Du98}
\be 
{\Omega}_2^{ij}(t)=(d-1+d_i+d_j)\Omega^{i\alpha}_1\Omega^{j\beta}_1
\partial_{t^\alpha}
\partial_{t^\beta} \mathbb{F}.
\ee

 Note that at this stage we are working under the hypothesis of Theorem \ref{main thm}. 
 
 Let us consider the coordinates \eqref{change coord} on $M\backslash \{t^1=0\}$. Then  the nonzero entries of the Jacobian matrix  are 
 \begin{align*}\label{Jvalue}
 \frac{\partial s^{i}}{\partial t^{1}}&=\frac{d_1-2d_i}{d_1}  t^i (t^1)^{\frac{-2d_i}{d_1}},\ \ 
\frac{\partial s^{r}}{\partial t^{1}}= (\frac{-2 -d_1}{2 d_1}) \sum_2^{r-1} t^{i} t^{r-i+1} (t^1)^{\frac{-2 }{d_1}-2}-\frac{2}{d_1} t^r (t^1)^{\frac{-2}{d_1}-1}, \\ \nonumber
\frac{\partial s^{i}}{\partial t^{i}}&= (t^1)^{\frac{d_1-2d_i}{d_1}}, \ \
\frac{\partial s^{r}}{\partial t^{i}}= t^{r-i+1} (t^1)^{\frac{-2 }{d_1}-1},\ \
\frac{\partial s^{r}}{\partial t^{r}}=  (t^1)^{\frac{-2 }{d_1}}.
 \end{align*}
 
\bp\label{new flat coord}
Consider the conjugate QFPM $\widetilde T=(\O_2,\widetilde \O_1)$. Then  $\widetilde\tau=s^1$, $\widetilde \Omega_1^{ij}(s)=\delta^{i+j}_{r+1}$, $\widetilde{e}=\partial_{s^r}$ and $\widetilde E=\widetilde d_i s^i\partial_{s^i}$ where the numbers $\widetilde d_i$ are given in \eqref{degrees inv}.
\ep 

\begin{proof}
 Using the duality between the degrees outlined in Lemma \ref{degrees},  we calculate the entries  $\widetilde\Omega_1^{ij}(s)$ as follows. 
\begin{enumerate}
    \item[I)] For $i=1$   \[\widetilde{\Omega}_1^{1j}(s)= - \frac{\partial s^{j}}{\partial t^{\alpha}} \widetilde{\Omega}_1^{1\alpha}= - \frac{\partial s^{r}}{\partial t^{r}}\widetilde{\Omega}_1^{1r} =- \frac{\partial s^{r}}{\partial t^{r}}(-(t^1)^{\frac{2}{d_1}}) \delta^{1r}=\delta^{1}_r.\]

\item[II)] For  $1 < i < r $ and $1< j< r $
\begin{align*}
\widetilde{\Omega}_1^{ij}(s)&=\frac{\partial s^{i}}{\partial t^{k}} \frac{\partial s^{j}}{\partial t^{k}} \widetilde{\Omega}_1^{k l}\\
 &=\frac{\partial s^{i}}{\partial t^{1}} \frac{\partial s^{j}}{\partial t^{1}} \widetilde{\Omega}_1^{11}+ \frac{\partial s^{i}}{\partial t^{i}} \frac{\partial s^{j}}{\partial t^{1}} \widetilde{\Omega}_1^{i1}+ \frac{\partial s^{i}}{\partial t^{1}} \frac{\partial s^{j}}{\partial t^{j}} \widetilde{\Omega}_1^{1j}+
 \frac{\partial s^{i}}{\partial t^{i}} \frac{\partial s^{j}}{\partial t^{j}} \widetilde{\Omega}_1^{ij}\\
 &= \frac{\partial s^{i}}{\partial t^{i}} \frac{\partial s^{j}}{\partial t^{j}} \widetilde{\Omega}_1^{ij} \delta^{i+j,r+1}\\
 &= (t^1)^{\frac{2d_1-2d_i-2d_{r-i+1}+2}{d_1}} \delta^{i+j,r+1}\\
 &= \delta^{i+j,r+1}.
 \end{align*}
 
 \item[III)] For $1< i < r$
\begin{align*}\widetilde{\Omega}_1^{i r}(s)&= (t^1)^{\frac{2}{d_1}} \frac{\partial s^{i}}{\partial t^{i}}\frac{\partial s^{r}}{\partial t^{r-i+1}}+ \left(-(t^1)^{\frac{2}{d_1}} \frac{\partial s^{i}}{\partial t^{1}}+\frac{-2 d_i}{d_1} t^i (t^1)^{\frac{2}{d_1}-1} \frac{\partial s^{i}}{\partial t^{i}}\right).\frac{\partial s^{r}}{\partial t^{r}}\\
&=(t^1)^{\frac{2}{d_1}} (t^1)^{\frac{d_1-2d_i}{d_1}}.t^{i} (t^1)^{\frac{-2 }{d_1}-1}+ \left(-\frac{d_1-2d_i}{d_1} (t^1)^{\frac{2}{d_1}} 
 t^i (t^1)^{\frac{-2d_i}{d_1}}+\frac{-2 d_i}{d_1} t^i (t^1)^{\frac{2}{d_1}-1} (t^1)^{\frac{d_1-2d_i}{d_1}}\right) (t^1)^{\frac{-2 }{d_1}}\\
&= (t^1)^{\frac{-2d_i }{d_1}}t^{i}
+ \left(-\frac{d_1-2d_i}{d_1} (t^1)^{\frac{-2d_i}{d_1}} 
 t^i +\frac{-2 d_i}{d_1}  (t^1)^{\frac{-2d_i}{d_1}}t^i \right)\\
 &= (t^1)^{\frac{-2d_i }{d_1}}t^{i}
- (t^1)^{\frac{-2d_i}{d_1}}t^{i} \\
&=0.
\end{align*}
\item[IV)] Finally, 
\begin{align*}\widetilde{\Omega}_1^{rr}(s)&= -(t^1)^{\frac{2}{d_1}} \frac{\partial s^{r}}{\partial t^{r}}.\frac{\partial s^{r}}{\partial t^{1}}+\sum_{i=2}^{r-1}  \left( (t^1)^{\frac{2}{d_1}} \frac{\partial s^{r}}{\partial t^{r-i+1}} -\frac{2 d_i}{d_1} t^i (t^1)^{\frac{2}{d_1}-1} \frac{\partial s^{r}}{\partial t^{r}}\right).\frac{\partial s^{r}}{\partial t^{i}}\\
&+\left( 
-(t^1)^{\frac{2}{d_1}} \frac{\partial s^{r}}{\partial t^{1}}+\sum_{i=2}^{r-1}-\frac{2 d_i}{d_1} t^i (t^1)^{\frac{2}{d_1}-1} \frac{\partial s^{r}}{\partial t^{i}} + \frac{-4}{d_1} t^r (t^1)^{\frac{2}{d_1}-1} \frac{\partial s^{r}}{\partial t^{r}}\right).\frac{\partial s^{r}}{\partial t^{r}}\\
&=(\frac{2}{d_1}+1) \sum_2^{r-1} t^i t^{r-i+1} (t^1)^{\frac{-2}{d_1}-2} +\frac{4}{d_1} t^r (t^1)^{\frac{-2}{d_1}-1}+\sum_2^{r-1} t^i t^{r-i+1} (t^1)^{\frac{-2}{d_1}-2}\\
&-\sum_2^{r-1} \frac{2d_i}{d_1}t^i t^{r-i+1} (t^1)^{\frac{-2}{d_1}-2}-\sum_2^{r-1}\frac{2d_{r-i+1}}{d_1} t^i t^{r-i+1} (t^1)^{\frac{-2}{d_1}-2}-\frac{4}{d_1} t^r (t^1)^{\frac{-2}{d_1}-1}\\
&=\sum_2^{r-1} \left( \frac{2}{d_1}+2 -\frac{2d_i}{d_1}-\frac{2d_{r-i+1}}{d_1} \right) t^i t^{r-i+1} (t^1)^{\frac{-2}{d_1}-2}\\
&=0.
\end{align*}
\end{enumerate}
  It is straightforward to show that  $\widetilde{e}=\partial_{s^r}$. The  vector field $\widetilde E=\Omega_2^{1j}(s)\partial_{s^j}$ while  
  \begin{align} \nonumber
{\Omega}_2^{1j}(s)&=\begin{pmatrix}  d_1 t^1& -d_1 t^1 \frac{\partial s^2}{\partial t^1} -d_2 t^2 \frac{\partial s^2}{\partial t^2}& &-d_1 t^1 \frac{\partial s^3}{\partial t^1} -d_3 t^3 \frac{\partial s^3}{\partial t^3}&&\cdots & -d_1 t^1\frac{\partial s^r}{\partial t^1} -d_2 t^2\frac{\partial s^r}{\partial t^2}+\cdots -t^r \frac{\partial s^r}{\partial t^r}
 \end{pmatrix} \\ \nonumber
 &= \begin{pmatrix} d_1 t^1& (d_2-d_1) t^2 (t^1)^{\frac{d_1-2d_2}{d_1}}&(d_3-d_1) t^3 (t^1)^{\frac{d_1-2d_3}{d_1}}&\cdots& \sum_{i=1}^{r} (-d_1 (\frac{-2-d_1}{2 d_1}) -d_i )t^i t^{r-i+1} (t^1)^{\frac{-2}{d_1}-1}
 \end{pmatrix} \\ \label{g(1j)}
 &= \begin{pmatrix} d_1 t^1& (d_2-d_1) t^2 (t^1)^{\frac{d_1-2d_2}{d_1}}&(d_3-d_1) t^3 (t^1)^{\frac{d_1-2d_3}{d_1}}&\cdots& \frac{1}{2}\sum_{i=1}^{r} t^i t^{r-i+1} (t^1)^{\frac{-2}{d_1}-1}\end{pmatrix}\\ \nonumber
 &= \begin{pmatrix}  - d_1 s^1 & \ \ \ (d_2-d_1) s^2  & \ \ \ \ \ \ \ \ \ \ \ (d_3-d_1) s^3 &  & \ \ \  \ \ \  \cdots & &\ \   s^r \ \ \  \ \ \ \end{pmatrix}.
 \end{align}
\end{proof}

We observe that the inverse transformation of  the inversion symmetry \eqref{Dob coord} is given by
\[ t^1=\frac{-1}{z^1},\ \ t^r=z^r + \frac{1}{2} \sum_2^{r-1} \frac{ z^i z^{r-i+1}}{z^1},\ \ t^k= \frac{-z^k}{z^1}, ~2\leq k\leq r.\]  
Thus, the potential \eqref{inv potential}   obtained from applying the inversion symmetry to  $\F(t)$ has the form \[\widetilde{\F}(z) = (z^1)^{2}    \F\left(\frac{-1}{z^1},\frac{-z^2}{z^1},\ldots,\frac{-z^{r-1}}{z^1},\frac{1}{2} \sum_1^r \frac{z^i z^{r-i+1}}{z^1}\right) +\frac{1}{2} z^r \sum_{1}^{r} z^i z^{r-i+1}.\]

\bl\label{potential in inv coord}
The potential $\widetilde\F(z)$ has the form
\be\label{F in 3 coord}
\widetilde{\F}(s) = (t^1)^{ \frac{-4}{d_1}}  \left(\F(t^1,\ldots,t^r) -\frac{1}{2} t^r \sum_{1}^{r} t^i t^{r-i+1}\right), z^i\leftrightarrow s^i.
\ee 
\el
\begin{proof}
We use the  identities 
 \[ t^1=- s^1= (s^1)^2 (\frac{-1}{s^1}),~ t^r= (s^1)^{\frac{2}{d_1}} \left(\frac{1}{2} \sum_1^r \frac{s^i s^{r-i+1}}{s^1} \right) ,~ t^i=(s^1)^{\frac{2 d_i}{d_1}} (\frac{-s^i}{s^1}), 1<i<r,\]
and  the quasihomogeneity of the potential $\F(t)$, i.e., 
 \be\label{quasi cond} (\frac{2}{d_1} E) \F(t)= \frac{2(3-d)}{d_1}  \F(t) =(\frac{4}{d_1}+2) \F(t).\ee
Then 
 \begin{align*}
&(t^1)^{ \frac{-4}{d_1}}  \left[\F(t^1,\ldots,t^r) -\frac{1}{2} t^r \sum_{1}^{r} t^i t^{r-i+1}\right]\\ \nonumber
&=(t^1)^{ \frac{-4}{d_1}}  \left[\F(t^1,\ldots,t^r)+\big( -t^1 (t^r)^2\big)-\frac{1}{2} t^r \sum_{2}^{r-1} t^i t^{r-i+1}\right]\\ \nonumber
    &=(s^1)^{ \frac{-4}{d_1}}  \left[\F \left( (s^1)^2 (\frac{-1}{s^1}),(s^1)^{\frac{2 d_2}{d_1}} (\frac{-s^2}{s^1}),\ldots,(s^1)^{\frac{2 d_{r-1}}{d_1}} (\frac{-s^{r-1}}{s^1}),(s^1)^{\frac{2}{d_1}} (\frac{1}{2} \sum_1^r \frac{s^i s^{r-i+1}}{s^1}) \right) +\big((s^r)^2 (s^1)^{\frac{4}{d_1}+1}\right. \big.\\ \nonumber
    & \left.\big. +s^r \sum_2^{r-1} s^i s^{r-i+1} (s^1)^{\frac{4}{d_1}} + s^1 \left( \frac{1}{2}\sum_{2}^{r-1} (s^1)^{\frac{2}{d_1}-1}  s^i s^{r-i+1}\right)^2\big)- \frac{1}{2}s^r (s^1)^\frac{4}{d_1}  \sum_{2}^{r-1} s^i s^{r-i+1} - s^1 \left( \frac{1}{2}\sum_{2}^{r-1} (s^1)^{\frac{2}{d_1}-1}  s^i s^{r-i+1}\right)^2\right]\\ \nonumber
        &=(s^1)^{\frac{-4}{d_1}} \left[(s^1)^{\frac{4}{d_1}+2} \F\left(\frac{-1}{s^1},- \frac{s^2}{s^1} ,-\frac{s^3}{s^1},\ldots, \frac{1}{2} \sum_{i=1}^{n} \frac{-s^i s^{n-i+1}}{s^1}  \right)   +(s^r)^2 (s^1)^{\frac{4}{d_1}+1}+ \frac{1}{2} s^r \sum_2^{r-1} s^i s^{r-i+1} (s^1)^{\frac{4}{d_1}} \right] \\ \nonumber
    &=(s^1)^{2}  \F\left(\frac{-1}{s^1},\frac{-s^2}{s^1},\ldots,\frac{-s^{r-1}}{s^1},\frac{1}{2} \sum_1^r \frac{s^i s^{r-i+1}}{s^1}\right)+ \frac{1}{2} s^r \sum_1^{r} s^i s^{r-i+1}
 \end{align*}
 which is the potential of the inversion symmetry by setting $s^i=z^i$.
 
 \end{proof}

Now we prove Theorem \ref{main thm} stated in the introduction. 
 
 \begin{proof}[Proof of Theorem \ref{main thm}] By Corollary \ref{reg bcz E} and Theorem \ref{dual Frob manif},  we use  the above notations and assume  $T=(\O_2,\O_1)$ is the associated QFPM. We need to show that the conjugate QFPM  $\widetilde T=(\O_2,\widetilde\O_1)$ equals the QFPM associated to the potential $\widetilde \F(s)$ given in \eqref{F in 3 coord}. This leads   to verifying that $\O_2 (s)$ equals the intersection form $\widehat{\O}_2(s)$ defined by  $\widetilde\F(s)$. It is straightforward to show that $\widetilde F(s)$ is a quasihomogenius function, i.e., $\widetilde E \widetilde F=(3-\widetilde d)\widetilde F$. Hence  \[
{\widehat\Omega}_2^{ij}(s):=(\widetilde d-1+\widetilde d_i+\widetilde d_j)\Omega^{i\alpha}_1\Omega^{j\beta}_1
\partial_{s^\alpha}
\partial_{s^\beta} \widetilde{\F}.
\]
After long calculations we find that $\widetilde{\Omega}_2^{ij}(s) =\widehat{\Omega}_2^{ij}(s)$. For examples,  
  we obtained the first row of  $ {\Omega}_2^{ij}(s)$  in \eqref{g(1j)} and  for even $r$ and $1< i,j< r$, we get by denoting $\partial_{t^i}\partial_{t^j} G(t)$  as $G_{i,j}$
\begin{align}\nonumber
{\Omega}_2^{ij} (s) &= \frac{\partial s^i}{\partial t^1} \frac{\partial s^j}{\partial t^1} \Omega_2^{1,1} +\frac{\partial s^i}{\partial t^i} \frac{\partial s^j}{\partial t^1} \Omega_2^{i,1} +
  \frac{\partial s^i}{\partial t^1} \frac{\partial s^j}{\partial t^j} \Omega_2^{1,j} +\frac{\partial s^i}{\partial t^i} \frac{\partial s^j}{\partial t^j} \Omega_2^{i,j} \\ \nonumber
  &=d_1(1-\frac{2d_i}{d_1})(1-\frac{2d_j}{d_1}) t^i t^{j} (t^1)^{1-\frac{2d_i}{d_1}-\frac{2d_j}{d_1}}+
  d_i (1-\frac{2d_j}{d_1}) t^i t^j (t^1)^{1-\frac{2d_i}{d_1}-\frac{2d_j}{d_1}}\\\nonumber
  &+ d_j (1-\frac{2d_i}{d_1}) t^i t^j (t^1)^{1-\frac{2d_i}{d_1}-\frac{2d_j}{d_1}}+
  (d-1+d_i+d_j) (t^1)^{2-\frac{2d_i}{d_1}-\frac{2d_j}{d_1}} (G_{r-i+1,n-j+1} + t^r \delta^{r,i+j})\\\nonumber
  &=(d_1-d_i-d_j) t^i t^j (t^1)^{1-\frac{2d_i}{d_1}-\frac{2d_j}{d_1}} + (-d_1+d_i+d_j) (t^1)^{2-\frac{2d_i}{d_1}-\frac{2d_j}{d_1}} \left(  G_{r-i+1,r-j+1} + t^r \delta^{r,i+j}\right)\\ \label{ex dual}
 &= (d_1-d_i-d_j) (t^1)^{1-\frac{2 d_i}{d_1}-\frac{2 d_j}{d_1}} \left( t^i t^j- t^1 G_{r-i+1,r-j+1} -t^1 t^r \delta^{r, i+j}  \right).
  \end{align}
  On the other hand
\begin{align}
    \frac{\partial^2 \widetilde{\F}}{\partial {s^{r-i+1}} \partial {s^{r-j+1}}}&= \left( t^r \delta_{r,i+j} (t^1)^{1-\frac{2}{d_1}-\frac{2d_{r-i+1}}{d_1}}+ G_{r-i+1,r-j+1} (t^1)^{-1-\frac{4}{d_1}+\frac{2d_{r-i+1}}{d_1}} \right) \left( -(s^1)^{\frac{2d_{r-j+1}}{d_1}-1} \right) \\ \nonumber
    &+ \left( t^i (t^1)^{1-\frac{2}{d_1}-\frac{2d_i}{d_1}} \right) \left( s^i (s^1)^{\frac{2}{d_1}-1} \right)\\ \nonumber
    &=\left( t^r \delta_{r,i+j} (t^1)^{2-\frac{2d_i}{d_1}-\frac{2d_j}{d_1}}+ G_{r-i+1,r-j+1} (t^1)^{-2-\frac{4}{d_1}+\frac{2d_{r-i+1}}{d_1}+\frac{2d_{r-j+1}}{d_1}} \right)- \left( t^i t^j (t^1)^{2-\frac{2d_i}{d_1}-\frac{2d_j}{d_1}} \right) \\ \nonumber
    &=(t^1)^{1-\frac{2d_i}{d_1}-\frac{2d_j}{d_1}} \left( t^r \delta^{r,i+j}  t^1 + G_{r-i+1,r-j+1} t^1-  t^i t^j  \right).
\end{align}
Therefore, 
\be \label{exx dual}
 \widehat{\Omega}_2^{ij}(s)= (d_i+d_j-d_1)(t^1)^{1-\frac{2d_i}{d_1}-\frac{2d_j}{d_1}} \left( t^r t^1 \delta^{r,i+j} + G_{r-i+1,r-j+1} t^1-  t^i t^j  \right) =\O_2^{ij}(s).
\ee
 \end{proof}
 
 \bx
Consider the following solution to WDVV equations
\be 
\F=\frac{t_1^3}{6}-\frac{1}{2} t_2^2 t_1+\frac{1}{2} t_2^2 t_3+\frac{1}{2} t_1 t_3^2.
\ee
It corresponds to a trivial  Frobenius manifold structure, i.e.,  Frobenius algebra structure does not depend on the point. Here the  charge $d=0$, the Euler vector field $E=\sum t_i \partial_{t_i}$ and identity vector field $e=\partial_{t_3}$. The intersection form is 
\[
\Omega_2(t) =\left(
\begin{array}{ccc}
 t_1 & t_2 & t_3 \\
 t_2 & t_3-t_1 & -t_2 \\
 t_3 & -t_2 & t_1 \\
\end{array}
\right)
\]
Setting 
\[
s_1=-t_1,\ \ s_2=\frac{t_2}{t_1},\ \  s_3=\frac{t_2^2}{2 t_1^3}+\frac{t_3}{t_1^2}
\]
the conjugate QFPM has $\widetilde \O_1^{ij}(s)=\delta^{i+j}_3$ and 
\[\O_2(s)=\left(
\begin{array}{ccc}
 -s_1 & 0 & s_3 \\
 0 & s_3+\frac{3 s_2^2}{2 s_1}+\frac{1}{s_1} & -\frac{s_2^3}{s_1^2}-\frac{2 s_2}{s_1^2} \\
 s_3 & -\frac{s_2^3}{s_1^2}-\frac{2 s_2}{s_1^2} & \frac{3 s_2^4}{4 s_1^3}+\frac{3 s_2^2}{s_1^3}-\frac{1}{s_1^3} \\
\end{array}
\right)\]
The potential of the  conjugate 
Frobenius manifold structure reads
\[\widetilde{\F}(s)=\frac{-1}{6 s_1} +\frac{s_2^2}{2 s_1}  +\frac{s_2^4}{8 s_1}+\frac{1}{2} s_2^2 s_3+\frac{1}{2} s_1 s_3^2.\]
One can check that this is  the same potential obtained by applying the inversion symmetry to $\F(t)$. Note that  ${\widetilde{E}}=-s_1 \partial_{s_1}+s_3 \partial_{s_3}$ and  $\widetilde{E} \widetilde{\F}= \widetilde{\F}$.
\ex

\section{The conjugate of a polynomial  Frobenius manifold} \label{first section}

In  this section, we  recall  the construction of Frobenius manifolds on the space of orbits of Coxeter groups given in \cite{DCG} and we apply the results of this article.

We fix an irreducible Coxeter group $\W$  of rank $r$. We consider the standard real reflection representation  $\psi: \W\to GL(V)$, where $V$ is  a complex vector space of dimension $r$.  Then the orbits space $M=V/\W$  is a variety whose coordinate ring is the ring of invariant polynomials $\C [V]^\W$. Using the Shephard-Todd-Chevalley theorem, the ring $\C [V]^\W$ is  generated by $r$ algebraically independent homogeneous polynomials. Moreover,   the degrees of a complete set of  generators are uniquely specified by the group \cite{Hum}.

We fix  a complete set of homogeneous  generators $u^1,u^2,\ldots,u^r$ for $\C [V]^\W$.  Let $\eta_i$ be the degree of $u^i$. Here, we have  \[2=\eta_1<\eta_2\leq\eta_3\leq  \ldots \leq\eta_{r-1}< \eta_r.\] It is known that $\eta_i+\eta_{r-i+1}=\eta_r+\eta_1$. Consider  the invariant bilinear form on $V$ under  the action of $\W$. Then it defines a contravariant flat metric $\Omega_2$ on $M$ and we let $u^1$ equals its quadratic form.   We fix the vector field $e:=\partial_{u^r}$. There is another  flat contravariant metric $\Omega_1:=\L_{e}\O_2$ on $M$, which was  initially studied by K. Saito (\cite{Saito}, \cite{Saito1}) and it is called the Saito flat metric. Then $T:=(\O_2,\O_1)$ is a FPM and Dubrovin  proved  the following theorem.  
\bt\cite{Du98}   $T=(\Omega_2,\O_1)$  is a regular QFPM of charge $\frac{\eta_r-2}{\eta_r}$ and leads to a polynomial Frobenius manifold structure on $M$, i.e., the corresponding potential is a polynomial function in the flat coordinates.
\et 

We observe that  the polynomial Frobenius structure defined by   $T$ has  $\tau=\frac{1}{\eta_r} u^1$, the Euler vector field    $E=\frac{1}{\eta_r}\sum_i\eta_i u^i \partial_{u^i}$, the identity vector field $e$ and  degrees  $\frac{\eta_i}{\eta_r}$. Note that $E$ is independent of the choice of generators but $e$ is defined up to a constant factor. Thus, changing the set of generators will lead to an equivalent  Frobenius manifold structure \cite{DCG}. The following theorem was conjectured by Dubrovin  and proved by C. Hertling. 

\bt  \cite{Hert} \label{polyFrob2}
Any  semisimple polynomial Frobenius manifold with  positive degrees  is isomorphic to  a polynomial Frobenius structure constructed on the   orbits space of the standard real reflection representation of a finite irreducible Coxeter group.
\et 
 
Clearly,  $T$ satisfies the hypotheses of Theorem \ref{dual Frob manif} and we have a conjugate regular QFPM $\widetilde T:=(\O_2, \L_{\widetilde e} \O_2)$, where  $\widetilde{e}= (\tau)^{\eta_r} e$.  Moreover,  from the work of K. Saito and his collaborators  (see also \cite{DCG}), we can fix  $u^1,\ldots, u^r$ to be flat with respect to $\Omega_1$ and the potential of the polynomial Frobenius manifold will  have the standard  form \eqref{norm potential}. In particular $\widetilde T$ is  the regular QFPM of the Frobenius manifold structure  obtained by applying inversion symmetry to  the polynomial Frobenius manifold on $M$. Considering Theorem \ref{polyFrob2}, we wonder what is the intrinsic description for the conjugate Frobenius manifold as this may  help in the classification  of Frobenius manifolds. 

In \cite{nonref}, we give a similar discussion for the $r$ Frobenius manifold structures constructed  in \cite{polyZuo} on the orbits space $M$ when $\W$ is of type $B_r$ or $D_r$.

 \section{Remarks}
 
  It is important to mention that the inversion symmetry of  the WDVV equation can be applied  to a solution  $\F(t)$ in the standard form \eqref{norm potential} under  more  general  quasihomogeneity condition than condition  \eqref{quasihomog} and without the  regularity condition \eqref{regcond} of the associated QFPM . In this case, if  the conjugate Frobenius manifold structure exists, we believe that it will be  equivalent to Frobenius manifold structure obtained by applying   the inversion symmetry,  we confirm this by Example \ref{example dim 2} and Example \ref{non regular but has new F}. 
  
  Note that Frobenius manifold structures which are invariant under inversion symmetry  were studied in \cite{morison}. We did not consider these cases as the charge will equal  1.
  
  It will be interesting to study the consequences of Theorem \ref{main thm}  on  the interpretation of the inversion symmetry in terms of the action of the Givental groups obtained in \cite{givental} and the  relation found in \cite{zang} between the principle hierarchies and tau functions of the two solutions to the WDVV equations related  by the inversion symmetry. We also believe that the findings in this article can be generalized to  the theory of bi-flat $F$-manifolds \cite{ArLor}. 
 
 It is known that the leading term of a certain class of compatible local Poisson  structures leads to a regular QFPM and thus to a Frobenius structure \cite{DZ}, \cite{Du98}.  Polynomial  Frobenius manifolds  obtained in \cite{mypaper1} are constructed by fixing  the regular nilpotent orbit in a simple Lie algebra and uses compatible local Poisson brackets obtained by Drinfeld-Sokolov reduction. In these cases, the  Poisson brackets form an exact Poisson pencil, and thus their central invariants are constants \cite{FalLor}. If the Lie algebra is simply-laced, then the central invariants  are equal \cite{DLZ} which means the Poisson structures  are consistent with  the principle hierarchy associated with the Frobenius manifold \cite{DZ}.   Fix one of these polynomial Frobenius structures and denote the associated  local Poisson brackets by  $\mathbb B_2$ and $\mathbb B_1$ (here $\mathbb B_2$ is the classical $W$-algebra). In the flat coordinates, these local Poisson brackets form an exact Poisson pencil under the identity vector field $e$, i.e.,  $\Lie_e\mathbb B_2=\mathbb B_1$ and $\Lie_e\mathbb B_1=0$. Let us denote the leading term of $\mathbb B_2$ by $B_2$ and  $\widetilde e$  is the vector field  associated with  the conjugate Frobenius manifold structure. We proved in this article that  $\Lie_{\widetilde e}^2 B_2=0$. Then   it is natural to ask if  $\widetilde e$ also leads to an exact Poisson pencil, i.e., $\Lie_{\widetilde e}^2\mathbb B_2=0$. Our calculations for the simple Lie algebra of type $A_3$, shows that this is not true.  
\vspace{0.1cm}

\noindent{\bf Acknowledgments.}
 The authors thank Paolo Lorenzoni for his time to read the first draft of this article and for his suggestion to generalize  the results under more general hypothesis. The authors also thank  anonymous reviewers whose comments/suggestions helped  clarify  and improve the article.  In particular, directing us to the coordinates free condition \eqref{new condition} and to the potential used in Example \ref{non regular but has new F}. 
\vspace{0.1cm}
 
\noindent{\bf Funding} This work was partially funded by the internal grant of Sultan Qaboos University \\ (IG/SCI/DOMS/19/08).

\vspace{0.1cm}

\noindent{\bf Data Availability} Non applicable. 
\section*{Declarations}

\noindent{\bf Conflict of interest} The authors have no relevant financial or non-financial interests to disclose.

\noindent Yassir Dinar \\
\noindent dinar@squ.edu.om \\

\noindent Zainab Al-Maamari\\
\noindent s100108@student.squ.edu.om\\

\noindent Depatment of Mathematics\\
\noindent College of Science\\
\noindent Sultan Qaboos University\\
\noindent Muscat, Oman

\end{document}